\documentclass[12pt]{article}

\usepackage{amssymb}
\usepackage{mathrsfs}

\usepackage{geometry}          

\usepackage{amsthm}
\usepackage{enumerate}

\usepackage{dsfont}
\usepackage{amsmath,amsthm,amscd,amssymb,bbm,graphicx}
\usepackage{floatflt,setspace, wrapfig}
\usepackage{showlabels}
\usepackage{graphics,multicol}
\usepackage{mathrsfs} 

\newtheorem{theorem}{Theorem}
\newtheorem{definition}[theorem]{Definition}
\newtheorem{lemma}[theorem]{Lemma}
\newtheorem{proposition}[theorem]{Proposition}
\newtheorem{corollary}[theorem]{Corollary}

\newcommand{\C}{\cal C}
\newcommand{\B}{\cal B}
\newcommand{\A}{\cal A}

\def\qed{\unskip\nobreak\hfill\penalty50\hskip 3pt\hbox{}\nobreak
\hfill\hbox{\vrule width 4 pt height 10 pt}}

\title{The Almost Sure Invariance Principle for Beta-Mixing measures}
\author{Nicolai Haydn\thanks{Mathematics Department, USC,
Los Angeles, 90089-1113. E-mail: $<$nhaydn@math.usc.edu$>$.}}

\begin{document}
\maketitle

\begin{abstract}
The theorem of Shannon-McMillan-Breiman states that for every generating
partition on an ergodic system of finite entropy
the exponential decay rate of the measure of cylinder sets 
equals the metric entropy almost everywhere. In addition the measure of $n$-cylinders 
is in various settings known to be lognormally distributed in the limit. 
In this paper the logarithm of the measure of $n$-cylinder, the information function,
satisfies the almost sure invariance principle in the case in which the measure
is $\beta$-mixing. We get a similar result for the recurrence time.
 Previous results are due to Philipp and Stout who deduced
the ASIP when the measure is strong mixing and satisfies an $\mathscr{L}^1$-type Gibbs
condition. 
\end{abstract}

\section{Introduction}

Let $\mu$ be a $T$-invariant probability measure on a space $\Omega$ on 
which the map $T$ acts measurably. For a  measurable partition $\cal A$
one forms the $n$th join ${\cal A}^n=\bigvee_{j=0}^{n-1}T^{-j}{\cal A}$
which forms a finer partition of $\Omega$. (The atoms of ${\cal A}^n$ are
traditionally called {\em $n$-cylinders}.)  For $x\in\Omega$
we denote by $A_n(x)\in{\cal A}^n$ the $n$-cylinder which contains $x$.
The Theorem of Shannon-McMillan-Breiman (see e.g.~\cite{Man,Petersen}) then states that for 
$\mu$-almost every $x$ in $\Omega$ the limit
$$
\lim_{n\rightarrow\infty}\frac{-\log\mu(A_n(x))}n
$$
exists and equals the metric entropy $h(\mu)$ provided the entropy is finite in the case
of a countable infinite partition. The convergence is uniform only in degenerate cases
(see~\cite{CGS} for an example).
 This theorem was proved for finite partitions in increasing degrees
of generality in the years 1948 to 1960 first for finite partitions and then for countably infinite partitions.
For a setting on metric spaces and with Bowen balls instead of cylinders, Brin and Katok~\cite{BK}
proved a similar almost sure limiting result.

Related to the SMB theorem are recurrence and waiting times for which limiting result
were proven by Ornstein and Weiss~\cite{OW,OW2} and Nobel and Wyner~\cite{NW} respectively. 
Here we are interested in a more detailed description of the limiting distribution of the
 information function $I_n(x)=-\log\mu(A_n(x))$ around its mean value.
These properties are of interest
when evaluating the efficiency of compression algorithms in information theory.

In 1962 Ibragimov~\cite{Ibr} proved the Central Limit Theorem for SMB for 
measures that are strongly mixing (in Rosenblatt's sense~\cite{R1}) and satisfy
an $\mathscr{L}^1$-type Gibbs condition, that is, he proved that $I_n$ is 
in the limit lognormally distributed. Various improvements followed although
most of them following Ibragimov's arguments or assume that the measure is Gibbs. For 
instance, Collet, Galves and Schmitt~\cite{CGS} proved that $I_n$ is lognormally 
distributed in the limit for $\psi$-mixing Gibbs measures, Paccot~\cite{Pac} for 
interval maps with suitable topological covering properties
 For other results see for instance~\cite{Bro,Pen,Gor,Liv,ch}. For Gibbs measures on non-uniformly 
expanding systems such results have been 
obtained in~\cite{FHV,BV}. For $(\psi,f)$-mixing measures a CLT was proven in~\cite{HV3},
for rational maps with critical points in the Julia set in~\cite{HV2}
and for $\beta$-mixing maps in~\cite{Hay}. This latter result does not require Ibragimov's
$\mathscr{L}^1$-Gibbs condition, but in return asks for the somewhat stronger mixing property,
 that is $\beta$-mixing instead of the strong mixing property.

In the setting of Ibragimov, Philipp and Stout~\cite{PS} then proved the 
almost sure invariance principle under similar conditions although with 
faster decay and better rates of approximability of the conditional entropy 
function. Kontoyiannis~\cite{Kon} then used this result to prove the almost sure invariance principle,
CLT and the law of the iterated logarithm LIL for recurrence and waiting times, 
strengthening the result of Nobel and Wyner~\cite{NW}.
Also Han~\cite{Han} proved the ASIP for SMB in the case of exponentially $\psi$-mixing systems
following Philipp and Stout. In the present paper we prove the ASIP for SMB 
for measures that are $\beta$-mixing. Here we don't require the $\mathscr{L}^1$-Gibbs
property of Ibragimov and Philipp and Stout. Also we allow for countably infinite partitions. 
These two aspects are the novelties of the present paper.
 
In Section~\ref{main.results} we define mixing conditions and state the main theorem.
In Section~\ref{variance&higher} we show existence of the variance and give the 
rate of convergence. We also obtain estimates on the growthrate of the higher order 
moments of the centred information function. These estimates are important in 
Section~\ref{proof.ASIP} where we proof the ASIP following the road laid out in~\cite{PS}.

\section{Main results}\label{main.results}
Let $T$ be a map on a space $\Omega$ and $\mu$ a
probability measure on $\Omega$. Moreover let $\cal A$ be a (possibly infinite)
 measurable partition of $\Omega$ and denote by
${\cal A}^n=\bigvee_{j=0}^{n-1}T^{-j}{\cal A}$
its {\em $n$-th join} which also is a measurable partition of $\Omega$ for
every $n\geq1$. The atoms of ${\cal A}^n$ are called {\em $n$-cylinders}.
Let us put ${\cal A}^*=\bigcup_{n=1}^\infty{\cal A}^n$ for the collection of
all cylinders in $\Omega$ and put $|A|$ for the length of a
cylinder $A\in{\cal A}^*$, i.e.\ $|A|=n$ if $A\in{\cal A}^n$.

We shall assume that $\cal A$ is generating, i.e.\ that the atoms of
${\cal A}^\infty$ are single points in $\Omega$.

\subsection{Mixing}

The main assumption in the results described here is on the mixing property of
 the invariant measure. Here we use the following:
 
\begin{definition}\label{definition}
We say the invariant probability measure $\mu$ is {\em $\beta$-mixing}\footnote{
 In~\cite{Hay} we used the term `uniform strong mixing' for what is commonly 
called $\beta$-mixing. Here we adhere to the standard terminology.}
if there exists a decreasing function $\psi: \mathbb{N}\rightarrow\mathbb{R}^+$ 
which satisfies $\beta(\Delta)\rightarrow0$ as $\Delta\rightarrow\infty$ so that
$$
\sum_{(B,C)\in {\cal A}^n\times T^{-\Delta-n}{\cal A}^m}\left|\mu(B\cap C)-\mu(B)\mu(C)\right|\le\beta(\Delta)
$$
for every every $n,m,\Delta>0$.
\end{definition}

\noindent {\bf Other kinds of mixing:}\\
For comparison purposes we list here some other kinds of mixing which are commonly used
in dynamics. Below $U$ is always in the $\sigma$-algebra
 generated by ${\cal A}^n$ and $V$ lies in the $\sigma$-algebra 
generated by ${\cal A}^*$ (see also~\cite{Dou}).
The limiting behaviour described is as the length of the `gap' $\Delta$ tends to infinity:
\begin{enumerate}
\item  {\em $\psi$-mixing}: 
$\displaystyle
\sup_n\sup_{U,V}\left|\frac{\mu(U\cap T^{-\Delta-n}V)}{\mu(U)\mu(V)}-1\right|\rightarrow0.
$
\item  {\em Left $\phi$-mixing}:
$\displaystyle
\sup_n\sup_{U,V}\left|\frac{\mu(U\cap T^{-\Delta-n}V)}{\mu(U)}-\mu(V)\right|\rightarrow0.
$
\item {\em Strong mixing} \cite{R1,Ibr} (also called $\alpha$-mixing):
$\displaystyle
\sup_n\sup_{U,V}\left|\mu(U\cap T^{-\Delta-n}V)-\mu(U)\mu(V)\right|\rightarrow0$.
\item   {\em Uniform mixing}~\cite{R1,R2}:
$\displaystyle
\sup_n\sup_{U,V}\left|\frac1k\sum_{j=1}^k\mu(U\cap T^{-n-j}V)-\mu(U)\mu(V)\right|\rightarrow0
$ as $k\rightarrow\infty$.

\end{enumerate}
One can also have {\em right $\phi$-mixing} when
$
\sup_n\sup_{U,V}\left|\frac{\mu(U\cap T^{-\Delta-n}V)}{\mu(V)}-\mu(U)\right|\rightarrow0
$
as $\Delta\rightarrow\infty$.
Clearly $\psi$-mixing is the strongest mixing property and  implies the other kinds of mixing.
The next strongest is $\phi$-mixing, then comes strong mixing and uniform 
mixing is the weakest. The $\beta$ mixing property is stronger that the 
strong mixing property but is implied by the $\phi$-mixing property. 

One says $\mu$ has the weak Bernoulli property (with respect to the 
partition $\cal A$) if for every $\varepsilon>0$ there exists an $N(\varepsilon)$ so that
$$
\sum_{B\in{\cal A}^n}\left|\mu(B\cap C)-\mu(B)\mu(C)\right|\le\varepsilon\mu(C)
$$
for every $C\in T^{-\Delta-n}{\cal A}^m$, $\Delta>N$ and $n,m\in\mathbb{N}$
(see e.g.~\cite{Petersen}). We see that the $\beta$-mixing property implies the weak
Bernoulli property. The rate $\beta$ determines
how fast the function $N(\varepsilon)$ grows as $\varepsilon$ goes to zero,
where to be precise $N(\varepsilon)=\beta^{-1}(\varepsilon)$.

\vspace{3mm}

\noindent For a partition $\cal A$ we have the ($n$-th)
{\em information function} $I_n(x)=-\log\mu(A_n(x))$, where $A_n(x)$ denotes the
unique $n$-cylinder that contains the point $x\in \Omega$, whose moments are
$$
K_w({\cal A}) = \sum_{A\in{\cal A}} \mu(A)|\log\mu(A)|^w=\mathbb{E}(I_n^w),
$$
$w\ge0$ not necessarily integer. (For $w=1$ one traditionally writes 
$H({\cal A})=K_1({\cal A})=\sum_{A\in{\cal A}}-\mu(A)\log\mu(A)$.) If $\cal A$ is
finite then $K_w({\cal A})<\infty$ for all $w$. For infinite partitions
the theorem of Shannon-McMillan-Breiman requires that $H({\cal A})$ be finite
to ensure finiteness of the entropy.
We will require finiteness of a larger than fifth moment  $K_w({\cal A})$ for some $w>5$ (not necessarily
integer).

\subsection{Results}
Our main result is the following theorem.

\begin{theorem}\label{ASIP}
Let  $\mu$ be a $\beta$-mixing invariant
probability measure on $\Omega$ with respect to a countably finite, measurable and 
generating partition $\cal A$ which satisfies $K_w({\cal A})<\infty$ for some $w>5$. 
Assume that $\beta$  decays at
 least polynomially with power $p>7+\frac{30}{w-5}$.

If the variance is positive, then the information function $I_n(x)=-\log\mu(A_n(x))$ satisfies the almost sure 
invariance principle for any error exponent $\delta\le\frac18$. 
That is, there exists a Brownian motion $B(n)$ such that 
$$
I_n=\sigma B(n)+\mathcal{O}(n^{\frac12-\delta})
$$
almost surely for any $\delta\le\frac18$. Moreover the variance $\sigma^2$ is given by
$$
 \sigma^2=\lim_{n\rightarrow\infty}\frac{K_2({\cal A}^n)-H^2({\cal A}^n)}n
 $$
 where the limit exists (and is strictly positive if the partition is infinite).
\end{theorem}

\noindent  Better bounds on $\delta$ are given in~\eqref{delta.bound}. 
The variance $\sigma^2$ is determined in Proposition~\ref{variance.convergence} 
and essentially only requires finiteness of the second moment $K_2({\cal A})$ although the 
rate of convergence uses that we have higher moments available. We obtain the following 
special cases using~\eqref{delta.bound}.

\begin{corollary} \label{corollary}
Suppose $\sigma>0$. Then:\\
(i) If $w=\infty$, e.g.\ if $\mathcal{A}$ is finite, and $\beta$ decays at  least polynomially with power $p>7$,
then $I_n=\sigma B(n)+\mathcal{O}(n^{\frac12-\delta})$
almost surely for any $\delta<\min\left(\frac12-\frac3{p+1},\frac13-\frac{10}{3p+30} \right)$.\\
(ii) If $\beta$ decays superpolynomially then $I_n$ satisfies the ASIP with exponent
$\delta<\frac13$ for any $w>5$.
\end{corollary}

\noindent Philipp and Stout~\cite{PS} Theorem~9.1 proved the ASIP under the condition that the measure
is strong mixing where $\phi(\Delta)=\mathcal{O}(\Delta^{-336})$ and requires that the 
$\mathscr{L}^1$-norms of the differences $f-f_n$ decay polynomially with power
$\ge 48$, where 
$f=\lim_{n\rightarrow\infty}f_n$ and 
$$
f_n(x)=\log\mathbb{P}(x_0|x_{-1}x_{-2}\dots x_{-n}).
$$
The ASIP holds then for any $\delta<\frac1{294}$. 

 If  $N(t)=\frac1{\sqrt{2\pi}}\int_{-\infty}^t e^{-s^2/2}\,ds$ denotes the normal distribution and
$h(\mu)$ the metric entropy of $\mu$ then we also have the CLT (\cite{Hay}):

\begin{theorem}\label{CLT}
Let  $\mu$ be a $\beta$-mixing 
probability measure on $\Omega$ with respect to a countably finite, measurable and 
generating partition $\cal A$ which satisfies $K_w({\cal A})<\infty$ for some $w>4$. 
Assume that $\beta$  decays at
 least polynomially with power $>6+\frac{20}{w-4}$.

If $\sigma>0$ then
$$
\mathbb{P}\left(\frac{I_n-nh}{\sigma\sqrt{n}}\le t\right)=N(t)+\mathcal{O}(n^{-\kappa})
$$
for all $t$ and all \\
(i) $\kappa<\frac1{10}-\frac35\frac{w}{(p+2)(w-2)+6}$ if $\beta $
 decays polynomially with power $p$,\\
(ii) $\kappa<\frac1{10}$ if $\beta$ decays super polynomially.
\end{theorem} 

\noindent 
The limiting result follows immediately from Theorem~\ref{ASIP} and the rate of 
convergence was obtained in~\cite{Hay} using Stein's method. There is an incomplete
argument in the variance and higher order estimates in~\cite{Hay} which are presented here 
in complete form (and also because we need higher than fourth moment)
and in fact here we obtain better lower bounds on the power $p$ than claimed in~\cite{Hay}.
 As a consequence, the error term
 for the variance has power $\frac14$ or better. All other estimates remain unchanged.
Another consequence of Theorem~\ref{ASIP} is  the Law of the Iterated Logarithm:

\begin{corollary}
Under the assumptions of Theorem~\ref{CLT}:
$$
\limsup_{n\rightarrow\infty}\frac{I_n(x)-nh}{\sigma\sqrt{2n\log\log n}}=1
$$
almost everywhere.
\end{corollary}

\noindent And similarly for the $\liminf$ where the limit  then equals $-1$ almost everywhere.
 Also in~\cite{Hay} we had proven the weak invariance principle 
WIP which then required to prove tightness and independence. It now follows directly from 
Theorem~\ref{ASIP} (although under slightly stronger assumptions).

\subsection{Examples}\label{examples}
{\bf (I) Bernoulli shift:} For the Bernoulli measure $\mu$ over the full shift space
$\Sigma=\mathbb{N}^\mathbb{Z}$ over the infinite alphabet $\mathbb{N}$
 generated by the weights $p_1,p_2,\dots$
($\sum_jp_j=1$), the entropy is then $h(\mu)=\sum_jp_j|\log p_j|$ and the variance is
$$
\sigma^2=\frac12\sum_{i,j}p_ip_j\log^2\frac{p_i}{p_j}
$$
assuming that $\sum_ip_i\log^2p_i<\infty$.
If moreover $\sum_ip_i|\log p_i|^{5+\epsilon}<\infty$ for some $\epsilon>0$ then 
we conclude the ASIP for $I_n(x)=-\log\mu(A_n(x))$ as $\beta$ decays exponentially fast.
\vspace{3mm}

\noindent
{\bf (II) Markov shift:} If $\mu$ is the Markov measure on $\Sigma=\mathbb{N}^\mathbb{Z}$ 
 generated by an infinite probability vector $\vec{p}=(p_1,p_2,\dots)$
($p_j>0$, $\sum_jp_j=1$) and an infinite stochastic matrix $P$ ($\vec{p}P=\vec{p}$, $P\mathbf{1}=\mathbf{1}$) then the entropy is $h(\mu)=\sum_{i,j}-p_iP_{ij}\log P_{ij}$~\cite{Wal} and
the variance~\cite{Rom,Yus,Hay} is
$$
\sigma^2=\frac12\sum_{ijk\ell}p_iP_{ij}p_kP_{k\ell}\log^2\frac{P_{ij}}{P_{k\ell}}
+4\sum_{k=2}^\infty\sum_{\vec{x}\in{\cal A}^k}\mu(\vec{x})\left(\log P_{x_1x_2}\log P_{x_{k-1}x_k}-h^2\right)
$$
where the terms in brackets on the RHS decay exponentially fast.
Then if $\sum_{i,j}p_iP_{i,j}|\log p_iP_{i,j}|^{5+\epsilon}<\infty$ for some $\epsilon>0$ then
$I_n$ satisfies the ASIP as $\beta$ decays exponentially fast ($p=\infty$).
Naturally we get the ASIP for any Markov measure over a finite alphabet.

\subsection{Recurrence time}

We denote by
$$
R_n(x)=\min\{k\ge1: T^kx\in A_n(x)\}
$$
 the $n$th {\em recurrence time} of $x$.
For a symbolic system where $T$ is the shift map on a symbolic space 
$\Sigma\subset \mathcal{A}^\mathbb{Z}$ the recurrence time
$R_n(\vec{x})=\min\{k\ge1: x_k\cdots x_{k+n}=x_0\cdots x_{n}\}$ of the point 
$\vec{x}(\dots,x_0,x_1,\dots)$ is the time is takes to see the first word of length $n$ again
as one goes to the `right'. Ornstein and Weiss~\cite{OW,OW2} showed that for 
ergodic measures $\lim_{n\to\infty}\frac1n|\log R_n(x)|=h(\mu)$ almost everywhere
improving on~\cite{WZ} where the convergence was shown to be in measure.
Collet,  Galves and Schmitt~\cite{CGS} proved the central limit theorem for
Gibbs measures which are exponentially $\psi$-mixing.
For finite alphabet processes Kontoyiannis~\cite{Kon} then proved the ASIP (for $\delta<\frac1{294}$)
 under the assumption that $\mu$
be $\alpha$-mixing with $\alpha$ decaying at least with power $336$ and that 
$\|f-f_n\|_{\mathscr{L}^1}$ decays with power $48$. Here we obtain the
following result which frees us from any condition on the sequence $\{f_n:n\}$.

\begin{theorem}
Assume $K_w(\mathcal{A})<\infty$ for some $w>5$, $\mu$ is $\beta$-mixing where
 $\beta$ decays at least least polynomially with power $p>7+\frac{30}{w-5}$
 and $\sum_n\|f-f_n\|_{\mathscr{L}^1}<\infty$.
 If $\sigma>0$ then 
 $$
 R_n=\sigma B(n)+\mathcal{O}(n^{\frac12-\delta})
 $$
 almost surely for all $\delta<\frac18$.
\end{theorem}

\noindent Better bounds for $\delta$ are as in~\eqref{delta.bound}. This theorem follows 
from Theorem~\ref{ASIP} and Corollary~1 of~\cite{Kon} where it was shown that 
$\log(R_n(x)\mu(A_n(x))=\log R_n-I_n=\mathcal{O}(n^\beta)$ for any $\beta>0$ provided
$\sum_n\|f-f_n\|_{\mathscr{L}^1}<\infty$. The classical case when $\mathcal A$ is a finite
alphabet requires that $\beta(\Delta)=\mathcal{O}(\Delta^{-p})$ for some $p>7$ and 
allows values $\delta<\min\left(\frac12-\frac3{p+1},\frac13-\frac{10}{3p+30} \right)$.

\section{Variance and higher moments}\label{variance&higher}

Before we prove the existence of the variance and bound the higher moments of the 
centred information function we shall summarise some known results which 
will be needed along the way.

\subsection{The information function}

\noindent   Denote by $A_n(x)$ the atom in ${\cal A}^n$ ($n=1,2,\dots$) which contains the point 
$x\in\Omega$. It was then shown in~\cite{GS} for $\psi$-mixing measures that 
$\sup_{A\in{\cal A}^n}\mu(A)$ decays exponentially fast as $n\to\infty$. For $\phi$-mixing
measures this was shown in~\cite{Aba} if $\phi(k)$ decays exponentially but is not 
necessarily true otherwise.
In~\cite{Hay} we then showed that for a $\beta$-mixing measure $\mu$ one has:\\
(i) $\sup_{A\in{\cal A}^n}\mu(A)=\mathcal{O}(n^{-p})$
if $\beta$ is polynomially decreasing with exponent $p>0$;\\
(ii) $\sup_{A\in{\cal A}^n}\mu(A)=\mathcal{O}(\theta^{\sqrt{n}})$ for some $\theta\in(0,1)$
if $\beta$ is exponentially decreasing.

The metric entropy $h$ for the invariant measure $\mu$ is
 $h=\lim_{n\rightarrow\infty}\frac1nH({\cal A}^n)$, where
  $\cal A$ is a generating partition of $\Omega$ (cf.\ \cite{Man}), provided $H({\cal A})<\infty$.
 For $w\ge1$ put $\eta_w(t)=t\log^w\frac1t$ ($\eta_w(0)=0$) and define
$$
K_w({\cal B}) = \sum_{B\in{\cal B}} \mu(B)|\log\mu(B)|^w
 = \sum_{B\in{\cal B}} \eta_w(\mu(B))
$$
for partitions $\cal B$.
Similarly one has the conditional quantity ($\cal C$  is a partition):
$$
K_w({\cal C}|{\cal B})
=\sum_{B\in{\cal B},C\in{\cal C}}\mu(B)\eta_w\left(\frac{\mu(B\cap C)}{\mu(B)}\right)
=\sum_{B,C}\mu(B\cap C)\left|\log\frac{\mu(B\cap C)}{\mu(B)}\right|^w.
$$

In the following we always assume that if $K_w({\cal A})<\infty$
then we also have $\mu(A)\le e^{-w}\;\forall\;A\in{\cal A}$. This can be achieved by passing to
a higher join. The assumption is convenient as it allows to use convexity arguments which
are implicit in some of the properties and estimates we use. We will need the following result.

\begin{lemma}\label{general.subadditivity} \cite{HV2}
For any two partitions ${\cal B}, {\cal C}$ for which $K_w({\cal B}),K_w({\cal C})<\infty$

\noindent (i) $K_w({\cal C}|{\cal B}) \leq K_w({\cal C})$,

\noindent (ii) $K_w({\cal B}\vee{\cal C})^{1/w}
\le K_w({\cal C}|{\cal B})^{1/w}+K_w({\cal B})^{1/w}$,

\noindent (iii) $K_w({\cal B}\vee{\cal C})^{1/w}
\le K_w({\cal C})^{1/w}+K_w({\cal B})^{1/w}.$
\end{lemma}

\noindent In~\cite{Hay} it was shown that as a consequence
\begin{equation}\label{K.estimate}
K_w({\cal A}^n)\le C_2 n^w.
\end{equation}
The variance of the information function $I_n$ is given by  
$\sigma_n^2=\sigma^2(I_n)= K_2({\cal A}^n)-H_n^2$ where $H_n=\mathbb{E}(I_n)=H({\cal A}^n)$.
For a partition ${\cal B}$ we write$J_{\cal B}$ for the centered information function given by
  $J_{\cal B}(B)=-\log\mu(B)-H({\cal B})$, $B\in{\cal B}$ (i.e. $\int J_{\cal B}\,d\mu=0$).
   Its variance is   $\sigma^2({\cal B})=\sum_{B\in{\cal B}}\mu(B)J_{\cal B}(B)^2$.
For two partitions $\cal B$ and $\cal C$ we put
 $$
 J_{{\cal C}|{\cal B}}(B\cap C)=\log\frac{\mu(B)}{\mu(B\cap C)}-H({\cal C}|{\cal B})
 $$
 for $(B,C)\in{\cal B}\times{\cal C}$. That is
 $J_{{\cal C}|{\cal B}}= J_{{\cal B}\vee{\cal C}}-J_{\cal B}$
 and $\sigma({\cal C}|{\cal B})=\sigma(J_{{\cal C}|{\cal B}})$. By~\cite{Hay}
\begin{equation}\label{variance.join}
\sigma({\cal B}\vee{\cal C})
\leq\sigma({\cal C}|{\cal B})+\sigma({\cal B}).
\end{equation}

\noindent As a consequence of Lemma~\ref{general.subadditivity}(i) one also has
$K_w({\cal B}\vee{\cal C}|{\cal B})=K_w({\cal C}|{\cal B})\le K_w({\cal C})$
which in particular implies 
$\sigma({\cal B}\vee{\cal C}|{\cal B})=\sigma({\cal C}|{\cal B})\le
\sqrt{K_2({\cal C})}$.
Let us put $\rho(B,C)=\mu(B\cap C)-\mu(B)\mu(C)$.
The following technical lemma is central to get the variance of $\mu$ and bounds
on the higher moments of $J_n=I_n-H_n$.

\begin{lemma}\label{rho.remainder.estimate}
Let $\mu$ be $\beta$-mixing and assume that $K_w({\cal A})<\infty$. 
Then for every
$\gamma>1$ and $a\in[1,w)$ there exists a constant $C_1$ so that 
$$
\sum_{B\in{\cal A}^m,C\in T^{-\Delta-m}{\cal A}^n}\hspace{-5mm}\mu(B\cap C)\left|\log\left(1+\frac{\rho(B,C)}{\mu(B)\mu(C)}\right)\right|^a
\le C_1\left(\beta(\Delta)(m+n)^{(1+a) \gamma}+(m+n)^{a \gamma-w(\gamma-1)}\right)
$$
for $\Delta<\min(n,m)$ and for all $n=1,2,\dots$.
\end{lemma}

\noindent We also have the following estimates for the approximations of $H(\mathcal{A}^n)$.

\begin{lemma}
Under the assumptions of Lemma~\ref{rho.remainder.estimate} the following applies:\\
(I) For every $\gamma>1$ there exists a constant $C_2$ so that for all $n$:
\begin{equation}\label{entropy.additivity}
\left|H({\cal A}^n\vee T^{-\Delta-n}{\cal A}^n)-2H({\cal A}^n)\right|
\le C_2\left(\beta(\Delta)n^{2 \gamma}+n^{\gamma-(\gamma-1)w}\right).
\end{equation}
(II) There exists a constant $C_3$ so that
\begin{equation}\label{entropy.approximation}
\left|\frac1mH({\cal A}^m)-h\right|\le C_3\frac1{m^\zeta}
\end{equation}
for all $m$,
where\\
(i) $\zeta\in(0,1-\frac{2w}{p(w-1)})$ if $\beta$ decays polynomially with power 
$p>\frac{2w}{w-1}$,\\
(ii) $\zeta\in(0,1)$ if $\beta$ decays faster than polynomially. 
\end{lemma}

\subsection{The variance}
In this  section we prove the existence of the variance as given in Theorem~\ref{ASIP} and moreover obtain
convergence rates.

\begin{proposition}\label{variance.convergence}
Let $\mu$ be $\beta$-mixing and assume that $K_w({\cal A})<\infty$ for some $w>2$. 
Let $\alpha\in(0,\frac12)$ and assume that $\beta$ is at least polynomially decaying with 
power $p\ge\frac{3\frac{w}\alpha-w-1}{w-2}$.
Then the limit
$$
\sigma^2=\lim_{n\rightarrow\infty}\frac1n\sigma^2({\cal A}^n)
$$
exists and is finite. Moreover there exists a constant $C_4$ so that
$$
\left|\sigma^2-\frac{\sigma^2({\cal A}^n)}n\right|\le \frac{C_4}{n^{\frac12-\alpha}}.
$$
If the partition $\cal A$ is infinite, then $\sigma$ is strictly positive.
\end{proposition}

\noindent {\bf Proof.}  Let us put ${\cal B}={\cal A}^n, {\cal C}=T^{-n-\Delta}{\cal A}^n$. 
 The gap $\Delta$ will be chosen to be $[n^\alpha]$ for some $\alpha\in(0,\frac12)$. 
We also assume here that $\beta$ decays polynomially with power $p$, that is 
$\beta(\Delta)=\mathcal{O}(\Delta^{-p})$.
Then by~\eqref{entropy.additivity}
$$
H({\cal B}\vee{\cal C})
=2H({\cal A}^n)+{\cal O}\left(\beta(\Delta)n^{2 \gamma}+n^{\gamma-(\gamma-1)w}\right)
=2H({\cal A}^n)+{\cal O}\left(n^{-\alpha p+2 \gamma}+n^{\gamma-(\gamma-1)w}\right)
$$
where $\gamma>1$ is arbitrary. The optimal value for $\gamma$ is $\gamma=\frac{\alpha p+w}{w+1}$ 
($\alpha p>1$) which yields the exponent  $\theta_1=\frac{-\alpha p(w-1)+2w}{w+1}$. That is 
$H({\cal B}\vee{\cal C})=2H({\cal A}^n)+{\cal O}(n^{\theta_1})$.
We get for the variance
\begin{eqnarray*}
\sigma^2({\cal B}\vee{\cal C})&=&\sum_{B\in{\cal B},C\in{\cal C}}
\mu(B\cap C)\left(\log\frac1{\mu(B\cap C)}
-H({\cal B}\vee{\cal C})\right)^2\\
&=&\sum_{B,C}
\mu(B\cap C)\left(J_{\cal B}(B)+J_{\cal C}(C)
+{\cal O}(n^{\theta_1})
-\log\left(1+\frac{\rho(B,C)}{\mu(B)\mu(C)}\right)\right)^2.
\end{eqnarray*}
By Minkowski's inequality:
$$
\left|\sigma({\cal B}\vee{\cal C})-\sqrt{ N({\cal B},{\cal C})} \right|
\le c_1n^{\theta_1}+\sqrt{F({\cal B},{\cal C})}
$$
($c_1>0$) where (by Lemma~\ref{rho.remainder.estimate} with $a=2$)
$$
F({\cal B},{\cal C})= \sum_{B\in{\cal B},C\in{\cal C}}
\mu(B\cap C)\log^2\left(1+\frac{\rho(B,C)}{\mu(B)\mu(C)}\right)
\le c_2\left(\beta(\Delta)n^{3 \gamma}+n^{2 \gamma-(\gamma-1)w}\right),
$$
for any $\gamma>1$ which when optimised yields the value $\gamma=\frac{\alpha p+w}{w+1}$.
Then with $\theta_2=\frac{-\alpha p(w-2)+3w}{w+1}$ we get $F({\cal B},{\cal C})=\mathcal{O}(n^{\theta_2})$.
The principal term is
\begin{eqnarray*}
N({\cal B},{\cal C})&=& \sum_{B\in{\cal B},C\in{\cal C}}
\mu(B\cap C)\left(J_{\cal B}(B)+J_{\cal C}(C)\right)^2\\
&=&\sum_{B,C}\mu(B\cap C)
\left(J_{\cal B}(B)^2+J_{\cal C}(C)^2\right)+2R({\cal B},{\cal C}) \\
&=& \sigma^2({\cal B})+ \sigma^2({\cal C})+ 2R({\cal B},{\cal C}).
\end{eqnarray*}
Since $J_{\cal B}$ and $J_{\cal C}$ have average zero the remainder term is
\begin{eqnarray*}
R({\cal B},{\cal C})&=&\sum_{B\in{\cal B},C\in{\cal C}}
\mu(B\cap C)J_{\cal B}(B)J_{\cal C}(C)\\
&=&\sum_{B,C}(\mu(B)\mu(C)+\rho(B,C))J_{\cal B}(B)J_{\cal C}(C)\\
&=&\sum_{B,C}\rho(B,C)J_{\cal B}(B)J_{\cal C}(C).
\end{eqnarray*}
In order to estimate $R$, put ${\cal L}=\{(B,C)\in{\cal B}\times{\cal C}: \mu(B\cap C)\ge2\mu(B)\mu(C)\}$
and write $R=R^++R^-$ where
\begin{eqnarray*}
R^+({\cal B},{\cal C})&=&\sum_{(B,C)\in{\cal L}}\rho(B,C)J_{\cal B}(B)J_{\cal C}(C),\\
R^-({\cal B},{\cal C})&=&\sum_{(B,C)\in{\cal L}^c}\rho(B,C)J_{\cal B}(B)J_{\cal C}(C).
\end{eqnarray*}
For $(B,C)\in{\cal L}$ one has $\rho(B,C)=\mu(B\cap C)-\mu(B)\mu(C)\ge\frac12\mu(B\cap C)$
and therefore, using H\"older's inequality twice ($\frac1r+\frac1s=1$ and $\frac1u+\frac1v=1$
such that $su, sv\le w$),
\begin{eqnarray*}
|R^+({\cal B},{\cal C})|
&\le&\sum_{(B,C)\in{\cal L}}\frac{\rho(B,C)}{\mu(B\cap C)}\left|J_{\cal B}(B)J_{\cal C}(C)\right|\mu(B\cap C)\\
&\le&\left(\sum_{(B,C)\in{\cal L}}\left(\frac{\rho(B,C)}{\mu(B\cap C)}\right)^r\mu(B\cap C)\right)^\frac1r
\left(\sum_{(B,C)\in{\cal L}}\left|J_{\cal B}(B)J_{\cal C}(C)\right|^{s}\mu(B\cap C)\right)^\frac1s\\
&\le&2^\frac{1}r\left(\sum_{B,C}|\rho(B,C)|\right)^\frac1r
\left(\left(\sum_{B}\left|J_{\cal B}(B)\right|^{su}\mu(B)\right)^\frac1u
\left(\sum_{C}\left|J_{\cal C}(C)\right|^{sv}\mu(C)\right)^\frac1v\right)^\frac1s\\
&\le&c_3\beta(\Delta)^\frac1rn^2
\end{eqnarray*}
since $\frac{\rho(B,C)}{\mu(B\cap C)}\le1\, \forall\, (B,C)\in{\cal L}$,
 the $\beta$-mixing property $\sum_{B,C}|\rho(B,C)|\le\beta(\Delta)$
 and where we used the a priori estimates
 $\sum_{B}\left|J_{\cal B}(B)\right|^{su}\mu(B)=\mathcal{O}(n^{us})$ (similarly for the sum over $C$).
We proceed similarly for the second part of the error term using the a priori estimate
$\sigma^2({\cal A}^n)\le c_4n^2$:
\begin{eqnarray*}
|R^-({\cal B},{\cal C})|
&\le&\sum_{(B,C)\in{\cal L}^c}\frac{|\rho(B,C)|}{\mu(B)\mu(C)}
\left|J_{\cal B}(B)J_{\cal C}(C)\right|\mu(B)\mu(C)\\
&\le&\left(\sum_{(B,C)\in{\cal L}^c}\left(\frac{|\rho(B,C)|}{\mu(B)\mu(C)}\right)^2\mu(B)\mu(C)\right)^\frac12
\left(\sum_{(B,C)\in{\cal L}^c}\left|J_{\cal B}(B)J_{\cal C}(C)\right|^{2}\mu(B)\mu(C)\right)^\frac12\\
&\le&\left(\sum_{B,C}|\rho(B,C)|\right)^\frac12
\sigma({\cal B})\sigma({\cal C})\\
&\le&c_4\beta(\Delta)^\frac1rn^2,
\end{eqnarray*}
where we used that $\mu(B\cap C)<2\mu(B)\mu(C)$ implies $|\rho(B,C)|\le\mu(B)\mu(C)$.
Hence with $s=\frac{w}2$ which is the largest possible value so that $\frac1r=\frac{w-2}w$ is the smallest 
possible and also $u=v=2$ say,  we obtain that 
$|R(\mathcal{B},\mathcal{C})|=\mathcal{O}(n^{2-\alpha p\frac{w-2}w})$
and therefore
\begin{equation}\label{sigma.addition}
\sigma ({\cal B}\vee{\cal C})
\leq\sqrt{\sigma^2({\cal C})+ \sigma^2({\cal B})+\mathcal{O}(n^{2-\alpha p\frac{w-2}w})}
+ \mathcal{O}(n^{\theta_2/2}).
\end{equation}
To fill the gap of length $\Delta$ estimates~(\ref{variance.join}) and~(\ref{K.estimate}) yield
$$
|\sigma({\A}^{2n+\Delta})- \sigma({\B}\vee{\C})|
\le\sigma(T^{-n}{\A}^\Delta|{\B}\vee{\C})
\le \sqrt{K_2(T^{-n}{\A}^{\Delta})}
=\sqrt{K_2({\A}^{\Delta})}
\le c_6 \Delta
=\mathcal{O}(n^\alpha).
$$
as $\Delta=[n^\alpha]$. We want to demand that 
$n^{1-\alpha p\frac{w-2}{2w}},n^{\theta_2/2}= \mathcal{O}(n^\alpha)$ which is achieved if
 $p\ge p_1=\frac{2w}{w-2}\frac{1-\alpha}\alpha$ and $p\ge p_2=\frac{3\frac{w}\alpha-w-1}{w-2}$ 
 respectively. Since $p_2>p_1$ we get the assumption $p\ge p_2$. Then, 
as $\sigma({\cal B})= \sigma({\C})=\sigma_n=\sigma({\cal A}^n)$,
one has
$$
\sigma_{2n+[n^\alpha]}= \sqrt{2\sigma^2_n+\mathcal{O}(n^{2\alpha})}+\mathcal{O}(n^\alpha),
$$
Since $\alpha<\frac12$ one has $\sigma^2_k\le c_{10}k$ for all $k$ and some constant $c_{10}$. 

In order to get the 
rate of convergence let $n_0$ be given put recursively $n_{j+1}=2n_j+[n_j^\alpha]$ ($j=0,1,2,\dots$).
Then $2^jn_0\le n_j\le2^jn_0\prod_{i=0}^{j-1}\left(1+\frac12n_i^{\alpha-1}\right)$
where the product is bounded by
$$
1\le \prod_{i=0}^{j-1}\left(1+\frac12n_i^{\alpha-1}\right)
\le\prod_{i=0}^{j-1}\left(1+\frac{1}{n_0^{1-\alpha}2^{(1-\alpha)i+1}}\right)
\le\exp\frac{c_{10}}{n_0^{1-\alpha}}.
$$
On the other hand as 
$\sigma_{n_{j+1}}=\sqrt{2\sigma^2_{n_j}+\mathcal{O}(n_j^{2\alpha})}+\mathcal{O}(n_j^\alpha)$, which 
yields
$$
\sigma_{n_{j+1}}^2
=2\sigma^2_{n_j}+\mathcal{O}(n_j^{2\alpha})
+\mathcal{O}(n_j^\alpha)\sqrt{2\sigma^2_{n_j}+\mathcal{O}(n_j^{2\alpha})}
=\left(\sqrt2\,\sigma_{n_j}+\mathcal{O}(n_j^{\alpha})\right)^2
$$
and consequently 
$$
\sigma_{n_{j+1}}=\sqrt2\,\sigma_{n_j}+\mathcal{O}(n_j^{\alpha}).
$$
Thus 
$$
\sigma_{n_j}=2^\frac{j}2\sqrt{n_0}+\sum_{k=0}^{j-1}2^\frac{j-1-k}2\mathcal{O}(n_k^\alpha)
$$
and since $\sqrt{n_j}=2^\frac{j}2\sqrt{n_0}e^{\mathcal{O}(n_0^{\alpha-1})}$ we get (as $\alpha<\frac12$)
$$
\frac{\sigma_{n_j}}{\sqrt{n_j}}
=\frac{\sigma_{n_0}}{\sqrt{n_0}}e^{\mathcal{O}(n_0^{\alpha-1})}
+\sum_{k=0}^{j-1}\frac{2^\frac{j-1-k}22^{\alpha k}n_0^\alpha\mathcal{O}(1)}{j^\frac{j}2\sqrt{n_0}}
=\frac{\sigma_{n_0}}{\sqrt{n_0}}e^{\mathcal{O}(n_0^{\alpha-1})}
+\mathcal{O}(n_0^{\frac12-\alpha}).
$$
Hence
$$
\frac{\sigma^2_{n_j}}{n_j}
=\frac{\sigma^2_{n_0}}{n_0}+{\cal O}\left(\frac1{n_0^{\frac12-\alpha}}\right).
$$
Taking $\limsup$ as $j\rightarrow\infty$ and $n_0\rightarrow\infty$ shows that the limit 
$\sigma^2=\lim_n\frac{\sigma^2_n}n$ exists and satisfies moreover 
$\left|\sigma^2-\frac{\sigma^2_n}n\right|\le C_4n^{-({\frac12-\alpha})}$
for some $C_4$. 
 
Positivity of $\sigma$ in the case of an infinite partition was shown in~\cite{Hay}.
\qed

\vspace{3mm}

\noindent For finite partitions the measure has variance zero
if it is a Gibbs state for a potential which is a coboundary.

\subsection{Higher order moments}
\noindent We will need estimates on the  higher moments of $J_n$ which we denote by 
$$
M_\ell({\cal B})=\sum_{B\in{\cal B}}\mu(B)|J_{\cal B}(B)|^\ell.
$$
the $w$th (absolute) moment of the function $J_{\cal B}$.
By Minkowski's inequality (on $\mathscr{L}^\ell$ spaces)
$$
M_\ell^\frac1\ell({\cal B}\vee{\cal C})
=\sqrt[\ell]{\mu\left(\left|J_{{\cal C}|{\cal B}}+J_{\cal B}\right|^\ell\right)}
\le\sqrt[\ell]{\mu(|J_{{\cal C}|{\cal B}}|^\ell)}+\sqrt[\ell]{\mu(|J_{\cal B}|^\ell)}
=M_\ell^\frac1\ell({\cal C}|{\cal B})+M_\ell^\frac1\ell({\cal B}),
$$
where $M_\ell({\cal C}|{\cal B})=\sum_{B\in{\cal B},C\in{\cal C}}\mu(B\cap C)|J_{{\cal C}|{\cal B}}(B\cap C)|^\ell$
are the conditional moments.
It follows from Corollary~\ref{K.estimate} that the absolute moments for the
joins ${\cal A}^n$ have the rough a priori estimate $M_\ell({\cal A}^n)\le K_\ell({\cal A}^n)\le C_2n^\ell$.
The purpose of the 
next proposition is to reduce the power from $\ell$ to $\frac\ell2$.

\begin{proposition}\label{higher.moments}
Let $\mu$ be $\beta$-mixing and assume that $K_w({\cal A})<\infty$
 for some $w>4$. 
Also assume that $\beta$ decays at least polynomially with power $p$.

Let $\ell$ be an integer strictly smaller than $w$, then there exists a constant $C_5$ 
so that for all $q\le\ell$
$$
M_{q}({\cal A}^n)\le C_5n^\frac{q}2
$$
provided $p\ge\frac{w(\ell+2)-\ell}{w-\ell}$.
\end{proposition}

\noindent {\bf Proof.} The statement is true for $\ell=2$ by Proposition~\ref{variance.convergence}.
We prove the result by induction for integers $\ell<[w]$. Assume the estimate is true for $k\le\ell-1$ and
 we will now prove it for $\ell$.

With ${\cal B}={\cal A}^n$, ${\cal C}=T^{-\Delta-n}{\cal A}^n$ and the gap
 $\Delta=[\sqrt{n}]$ ($\alpha=\frac12$ in the previous notation)  we get
by ~\eqref{entropy.additivity}
 $H({\cal B}\vee{\cal C})
 =H({\cal B}) + H({\cal C})  + {\cal O}(\beta(\Delta)n^{2\gamma}+n^{\gamma-(\gamma-1)w})$
 for $\gamma>1$ which as in the proof of Proposition~\ref{variance.convergence} optimised gives 
 $H({\cal B}\vee{\cal C})=H({\cal B}) + H({\cal C})  + {\cal O}(n^{\theta_1})$
  ($\theta_1=\frac{-p(w-1)+4w}{2(w+1)}$).
With Minkowsky's inequality
\begin{eqnarray*}
M_{\ell}^\frac1{\ell}({\cal B}\vee{\cal C})&=&\left(\sum_{B\in{\cal B},C\in{\cal C}}
\mu(B\cap C)\left|\log\frac1{\mu(B\cap C)}
-H({\cal B}\vee{\cal C})\right|^{\ell}\right)^\frac1{\ell}\\
&=&\hat{M}_{\ell}^\frac1{\ell}({\cal B},{\cal C})+{\cal O}(n^{\theta_1})
+F_{\ell}^\frac1{\ell}({\cal B},{\cal C}),
\end{eqnarray*}
where 
$$
\hat{M}_\ell({\cal B},{\cal C})= \sum_{B\in{\cal B},C\in{\cal C}}
\mu(B\cap C)\left|J_{\cal B}(B)+J_{\cal C}(C)\right|^{\ell}
$$
and by Lemma~\ref{rho.remainder.estimate} (with $a=\ell$) we get the error estimate
$$
F_{\ell}({\cal B},{\cal C})= \sum_{B\in{\cal B},C\in{\cal C}}
\mu(B\cap C)\left|\log\left(1+\frac{\rho(B,C)}{\mu(B)\mu(C)}\right)\right|^\ell
= {\cal O}\left(\beta(\Delta)n^{(\ell+1) \gamma}+n^{\ell \gamma-(\gamma-1)w}\right)
$$
where the value of $\gamma>1$ is optimised when $\gamma=\frac{p+2w}{2(w+1)}$ which 
then yields $F_{\ell}({\cal B},{\cal C})=\mathcal{O}(n^{\theta_\ell})$ where 
$\theta_\ell=\frac{-p(w-\ell)+2(\ell+1)w}{2(w+1)}$.
We want to achieve that $n^{\theta_1}=\mathcal{O}(n^\frac12)$ and, more generally,
$n^{\theta_\ell}=\mathcal{O}(n^\frac\ell2)$. That is $\theta_\ell\le\frac\ell2$ 
and this is satisfied if (as by assumption) $p\ge p_1=\frac{w(\ell+2)-\ell}{w-\ell}$. 
Hence $F(\mathcal{B},\mathcal{C})=\mathcal{O}(n^\frac\ell2)$
and moreover
\begin{equation}\label{first.approximation}
M_\ell(^\frac1\ell{\cal B}\vee{\cal C})
=\hat{M}_\ell^\frac1\ell({\cal B},{\cal C})+\mathcal{O}(n^\frac12).
\end{equation}
We further approximate 
$$
\hat{M}_\ell({\cal B},{\cal C})=M_\ell^\times({\cal B},{\cal C})+R_\ell({\cal B},{\cal C})
$$
where 
$$
M_\ell^\times({\cal B},{\cal C})
= \sum_{B\in{\cal B},C\in{\cal C}}\mu(B)\mu(C)\left|J_{\cal B}(B)+J_{\cal C}(C)\right|^{\ell}
$$
is the principal term and the error term $R_\ell=R_\ell^++R_\ell^-$ is split into a sum of two terms as 
in the proof of Proposition~\ref{variance.convergence}. 
Put ${\cal L}=\{(B,C)\in{\cal B}\times{\cal C}: \mu(B\cap C)\ge2\mu(B)\mu(C)\}$ and then
\begin{eqnarray*}
R_\ell^+({\cal B},{\cal C})&=&\sum_{(B,C)\in{\cal L}}\rho(B,C)\left|J_{\cal B}(B)+J_{\cal C}(C)\right|^{\ell},\\
R_\ell^-({\cal B},{\cal C})&=&\sum_{(B,C)\in{\cal L}^c}|\rho(B,C)|\left|J_{\cal B}(B)+J_{\cal C}(C)\right|^{\ell}.
\end{eqnarray*}
We now estimate the two terms separately as follows:\\
{\bf (I)}
For $(B,C)\in{\cal L}$ one has $\rho(B,C)=\mu(B\cap C)-\mu(B)\mu(C)\ge\frac12\mu(B\cap C)$
and therefore, by H\"older's inequality ($\frac1s+\frac1r=1$),
\begin{eqnarray*}
R_\ell^+({\cal B},{\cal C})
&=&\sum_{(B,C)\in{\cal L}}\frac{\rho(B,C)}{\mu(B\cap C)}\left|J_{\cal B}(B)+J_{\cal C}(C)\right|^{\ell}\mu(B\cap C)\\
&\le&\left(\sum_{(B,C)\in{\cal L}}\left(\frac{\rho(B,C)}{\mu(B\cap C)}\right)^r\mu(B\cap C)\right)^\frac1r
\left(\sum_{(B,C)\in{\cal L}}\left|J_{\cal B}(B)+J_{\cal C}(C)\right|^{s\ell}\mu(B\cap C)\right)^\frac1s\\
&\le&2^\frac{1}r\left(\sum_{B,C}|\rho(B,C)|\right)^\frac1r
\hat{M}_{s\ell}({\cal B},{\cal C})^\frac1s\\
&\le&2^\frac{1}r\beta(\Delta)^\frac1r\hat{M}_{s\ell}({\cal B},{\cal C})^\frac1s
\end{eqnarray*}
($s\ell\le w$) since $\frac{\rho(B,C)}{\mu(B\cap C)}\le1\, \forall \,(B,C)\in{\cal L}$ and 
where we used the $\beta$-mixing property $\sum_{B,C}|\rho(B,C)|\le\beta(\Delta)$.
We now use the a priori estimate 
$M_{q}({\cal A}^n)\le c_2n^q$ to 
bound $\hat{M}_{s\ell}({\cal B},{\cal C})$. Using Minkowsky's inequality we get the rough a priori bound
for $q\le w$:
\begin{eqnarray*}
\hat{M}_{q}^\frac1q({\cal B},{\cal C})&=& \left(\sum_{B\in{\cal B},C\in{\cal C}}
\mu(B\cap C)\left|J_{\cal B}(B)+J_{\cal C}(C)\right|^q\right)^\frac1q\\
&\le& \left(\sum_{B\in{\cal B}}\mu(B)\left|J_{\cal B}(B)\right|^q\right)^\frac1q
+ \left(\sum_{C\in{\cal C}}\mu(C)\left|J_{\cal C}(C)\right|^q\right)^\frac1q\\
&=&2M_q^\frac1q({\cal A}^n)\\
&\le&2c_2n
\end{eqnarray*}
and consequently $R_\ell^+({\cal B},{\cal C})=\mathcal{O}(\beta(\Delta)^\frac1rn^\ell)$.\\
{\bf (II)} We proceed similarly for the second part of the error term ($\frac1s+\frac1r=1$):
\begin{eqnarray*}
|R_\ell^-({\cal B},{\cal C})|
&\le&\sum_{(B,C)\in{\cal L}^c}\frac{|\rho(B,C)|}{\mu(B)\mu(C)}\left|J_{\cal B}(B)+J_{\cal C}(C)\right|^{\ell}\mu(B)\mu(C)\\
&\le&\left(\sum_{(B,C)\in{\cal L}^c}\left(\frac{|\rho(B,C)|}{\mu(B)\mu(C)}\right)^r\mu(B)\mu(C)\right)^\frac1r
\left(\sum_{(B,C)\in{\cal L}^c}\left|J_{\cal B}(B)+J_{\cal C}(C)\right|^{s\ell}\mu(B)\mu(C)\right)^\frac1s\\
&\le&\left(\sum_{B,C}|\rho(B,C)|\right)^\frac1r
M_{s\ell}^\times({\cal B},{\cal C})^\frac1s\\
&\le&\beta(\Delta)^\frac1rM_{s\ell}^\times({\cal B},{\cal C})^\frac1s,
\end{eqnarray*}
where we used that $\mu(B\cap C)<2\mu(B)\mu(C)$ implies $|\rho(B,C)|\le\mu(B)\mu(C)$.
Using the a priori estimate $M_{q}({\cal A}^n)\le c_2n^{q}$, we obtain by Minkowsky for $q\le w$:
\begin{eqnarray*}
(M_{q}^\times({\cal B},{\cal C}))^\frac1q
&\le& \left(\sum_{B\in{\cal B}}\mu(B)\left|J_{\cal B}(B)\right|^q\right)^\frac1q
+ \left(\sum_{C\in{\cal C}}\mu(C)\left|J_{\cal C}(C)\right|^q\right)^\frac1q\\
&=&2M_q^\frac1q({\cal A}^n)\\
&\le&2c_2n
\end{eqnarray*}
and consequently $R_\ell^-({\cal B},{\cal C})=\mathcal{O}(\beta(\Delta)^\frac1rn^\ell)$.

\vspace{3mm}

\noindent The two parts~(I) and (II) combined yield
$$
|R_\ell({\cal B}\vee{\cal C})|\le c_3 n^{\ell-\frac{p}{2r}}
$$
where we choose $s$ and $r$ such that $s\ell\le w$ and $\ell -\frac{p}{2r}\le\frac\ell2$.
For $s=\frac{w}\ell$ (largest possible) and $r=\frac{w}{w-\ell}$ this requires $p\ge p_2=\frac{w\ell}{w-\ell}$
which is satisfied by the assumption since $p_1>p_2$.
Then
$$
M_{\ell}({\cal B}\vee{\cal C})=M_{\ell}^\times({\cal B},{\cal C})+\mathcal{O}(n^\frac\ell2)
$$
and now look more closely at the principal term $M_{\ell}^\times$. Using the induction hypothesis
 $M_k({\cal A}^n)=\mathcal{O}(n^\frac{k}2)$ for $k\le \ell-1$ we obtain for $\ell$ integer
$$
M_{\ell}^\times({\cal B},{\cal C})
\le M_\ell({\cal B})+M_\ell({\cal C})
+\sum_{k=1}^{\ell-1}\left(\begin{array}{c}\ell\\k\end{array}\right)M_k({\cal B})M_{\ell-k}({\cal C})
\le 2M_\ell({\cal A}^n)+c_4n^\frac\ell2
$$
and consequently
$$
M_{\ell}^\frac1\ell({\cal B}\vee{\cal C})
=\left(2M_{\ell}({\cal A}^n)+\mathcal{O}(n^\frac\ell2)\right)^\frac1\ell+\mathcal{O}(n^\frac12).
$$
To fill in the gap of length $\Delta$ we use Lemma \ref{general.subadditivity}(iii),
 the estimate~(\ref{K.estimate}) on $K_\ell$ and the fact that $\Delta\sim \sqrt{n}$:
$$
\left|M_\ell^\frac1\ell({\A}^{2n+\Delta})-M_\ell^\frac1\ell({\B}\vee{\C})\right|\le
M_\ell^\frac1\ell({\A}^{2n+\Delta} |{\B}\vee{\C})
\le K_\ell^\frac1\ell({\A}^\Delta)
=\mathcal{O}(\Delta)=\mathcal{O}(n^\frac12).
$$
which implies 
$$
M_\ell^\frac1\ell({\A}^{2n+\Delta})
\le\left(2M_{\ell}({\cal A}^n)+c_5n^\frac\ell2\right)^\frac1\ell+c_6n^\frac12
$$
for constants $c_5, c_6$. Given $n_0$, put recursively $n_{j+1}=2n_j+[\sqrt{n_j}]$ ($\Delta_j=[\sqrt{n_j}]$ are the gaps), then for a constant $c_7$ large enough so that 
$(2+c_5/c_7)^\frac1\ell+c_6/c_7^\frac1\ell\le\sqrt2$ (which is possible as $\ell>2$) we obtain
$$
M_\ell({\cal A}^{n_j})\le c_7n_j^\frac\ell2
$$
for all $j$. Increasing the constant $c_7$ allows us to extend the estimate to all $n$ with a 
constant $C_5$.
This completes the inductions step.
If $\ell$ is the largest integer strictly smaller than $w$, then we can use H\"older's inequality to 
extend the estimate $M_q(\mathcal{A}^n)\le c_7n^\frac{q}2$ to arbitrary values of  $q\le \ell$.
\qed

\section{Proof of the ASIP (Theorem~\ref{ASIP})}\label{proof.ASIP}

Let $\alpha$ denote a number between $0$ and $1$ and $\ell<w$ an integer.
We decompose $J_n=I_n-H(\mathcal{A}^n)=\sum_{j=1}^{Q_n}(y_j+z_j)$
where $|z_j|=\Delta_j$ and $\Delta_j=[n_j^\alpha]$ is the length 
of the gaps where the length $|y_j|=n_j$ will be chosen 
to be $[\sqrt{j}] $.
Then $n=\sum_{j=1}^{Q_n}(n_j+\Delta_j)+r_n$ with remainder $r_n<n_{Q_n+1}+\Delta_{Q_n+1}$.
In particular $n\asymp \sum_{j=1}^{Q_n}\sqrt j \asymp Q_n^\frac32$ implies that $Q_n \asymp n^\frac32$.
We put $\hat{\mathcal A}^n=\bigvee_{j=1}^{Q_n}T^{-N_j}\mathcal{A}^{n_j}$
where we put $N_j=\sum_{i=1}^{j-1}(n_i+\Delta_i)$
and $N_1=0$. We have $N_j \asymp \sum_{i=1}^{j-1}\sqrt i \asymp j^\frac32$.

\begin{lemma} Let $\ell\ge5$ and $p\ge\frac{2(w+2)-\ell}{w-\ell}$. Then
$$
J_n=\sum_{j=1}^{Q_n}J_{n_j}\circ T^{N_j}+\mathcal{O}(n^{1-\delta_1})
$$
almost surely for any $\delta_1\le\min(\frac{p(w-2)-5w+1}{2p(w-2)+2w+2},\frac{1-\alpha}3)$.
\end{lemma}

\noindent {\bf Proof.} We proceed in three steps. First we cut `gaps' of lengths $\Delta_j$,
then we use the $\beta$-mixing property to separate the long blocks of lengths $n_j$
and in the last part we adjust the averaging term (entropy).\\
{\bf (I)} Since $J_{\mathcal{A}^n|\hat{\mathcal{A}}^n}
=J_{\mathcal{A}^n\vee\hat{\mathcal{A}}^n}-J_{\hat{\mathcal{A}}^n}
=J_{\mathcal{A}^n}-J_{\hat{\mathcal{A}}^n}$ (as $\mathcal{A}^n$ refines $\hat{\mathcal{A}}^n$)
we obtain
$$
\|J_{\mathcal{A}^n}-J_{\hat{\mathcal{A}}^n}\|_a^a
=\|J_{\mathcal{A}^n|\hat{\mathcal{A}}^n}\|_a^a
=M_a(\mathcal{A}^n|\hat{\mathcal{A}}^n)
\le\sum_{j=1}^{Q_n}M_a(\mathcal{A}^{\Delta_j})
$$
and using Proposition~\ref{higher.moments} for $1<a\le\ell$ (as $p\ge \frac{2(\ell+2)-\ell}{w-\ell}$ 
by assumption)
$$
\|J_{\mathcal{A}^n}-J_{\hat{\mathcal{A}}^n}\|_a^a
\le c_1\sum_{j=1}^{Q_n}n_j^{\frac{a\alpha}2}
\le c_1\sum_{j=1}^{Q_n}j^{\frac{a\alpha}4}
\le c_2 Q_n^{\frac{a\alpha}4+1}
=\mathcal{O}(n^{\frac{a\alpha}6+\frac23}).
$$ 
By Tchebycheff
$$
\mathbb{P}(|J_{\mathcal{A}^n}-J_{\hat{\mathcal{A}}^n}|>\epsilon_n)
\le\frac{\|J_{\mathcal{A}^n|\hat{\mathcal{A}}^n}\|_a^a}{\epsilon_n^a}
\le c_3\frac{n^{\frac{a\alpha}6+\frac23}}{n^{a(1-\delta)}}
$$
where we put $\epsilon_n=n^{1-\delta}$. This is summable if $a(1-\delta)-\frac{a\alpha}6-\frac23>1$
which is satisfied we we choose $a=\ell$ and as $\ell\ge5$ this is satisfied for any $\alpha<1$ and 
$\delta<\frac13$.
Thus by Borel-Cantelli
$$
J_{\mathcal{A}^n}=J_{\hat{\mathcal{A}}^n}+\mathcal{O}(n^{1-\delta})
$$
almost surely for any $\delta\le\frac13$.\\
{\bf (II)} Now let us put $D_k=\bigvee_{j=1}^{k-1}T^{-N_j}\mathcal{A}^{n_j}$, i.e.\
recursively we have
$D_{k+1}=D_k\vee T^{-N_k}\mathcal{A}^{n_k}$ and
also 
$D_{Q_n+1}=\hat{\mathcal{A}}^n$. 
Then by Lemma~\ref{rho.remainder.estimate} (with identification 
$\Delta=\Delta_{k-1}, n=N_\ell-\Delta_{k-1}, m=n_k, n+m\le N_{k+1}$ and not necessarily the 
same number $a$ as in part~(I))
\begin{eqnarray*}
\|I_{D_{k+1}}-I_{D_k}-I_{n_k}\circ T^{N_k}\|_a^a\hspace{-3cm}\\
&=&\sum_{D\in D_k, A\in T^{-N_k}\mathcal{A}^{n_k}}
\mu(D\cap A)\left|\log\left(\frac1{\mu(D\cap A)}-\log\frac1{\mu(A)}-\log\frac1{\mu(D)}\right)\right|^a\\
&=&\mathcal{O}\left(\beta(\Delta_{k-1})N_{k+1}^{(1+\alpha)\gamma}
+N_{k+1}^{\alpha\gamma-w(\gamma-1)}\right)
\end{eqnarray*}
for any $\gamma>1$ and $1<a<w$. As $\beta(\Delta)=\mathcal{O}( \Delta^{-p})$
\begin{eqnarray*}
\|I_{D_{k+1}}-I_{D_k}-I_{n_k}\circ T^{N_k}\|_a^a
&=&\mathcal{O}\left(n_{k-1}^{-\alpha p}N_k^{(1+a)\gamma}+N_k^{a\gamma-w(\gamma-1)}\right)\\
&=&\mathcal{O}\left(k^{-\frac{\alpha p}2}k^{\frac32\gamma(1+a)}+k^{\frac32(a\gamma-w(\gamma-1)}\right)
\end{eqnarray*}
which implies by Minkowski that
\begin{eqnarray*}
\|I_{\hat{\mathcal{A}}^n}-\sum_{j=1}^{Q_n}I_{n_j}\circ T^{N_j}\|_a
&=&\mathcal{O}(1)\sum_{k=1}^{Q_n}\left(k^{(-\frac{\alpha p}2
+\frac32\gamma(1+a))\frac1a}+k^{(\frac32(a\gamma-w(\gamma-1))\frac1a}\right)\\
&=&\mathcal{O}\left(n^{(\gamma(1+a)-\frac{\alpha p}3+\frac23)\frac1a}
+n^{((a\gamma-w(\gamma-1)+\frac23)\frac1a}\right)
\end{eqnarray*}
Thus
$$
\mathbb{P}\left(\left|I_{\hat{\mathcal{A}}^n}-\sum_{j=1}^{Q_n}I_{n_j}\circ T^{N_j}\right|>\epsilon_n\right)
\le\frac1{\epsilon_n^a}\left\|I_{\hat{\mathcal{A}}^n}-\sum_{j=1}^{Q_n}I_{n_j}\circ T^{N_j}\right\|_a^a
$$
is summable over $n\in\mathbb{N}$ if ($\gamma>1, 1<a<w$)
$$
\left\{\begin{array}{rcl}
\gamma(1+a)-\frac13\alpha p+\frac23+\delta&<&0\\
a\gamma-w(\gamma-1)+\frac23+\delta&<&0
\end{array}\right.
$$
where we chose $\epsilon_n=n^{1-\delta}$. The conditions on $\gamma$ and $a$ are satisfied
for any $\delta\le\frac13$ since $w>5$ and $p>4$.
Hence by Borel-Cantelli
$$
\left|I_{\hat{\mathcal{A}}^n}-\sum_{j=1}^{Q_n}I_{n_j}\circ T^{N_j}\right|
\le n^{1-\delta}
$$
for all $n$ large enough almost everywhere.\\
{\bf (III)} The entropies are estimated using Lemma~\ref{entropy.approximation} as follows:
$$
\left|H(\mathcal{A}^n)-\sum_{j=1}^{Q_n}H(\mathcal{A}^{n_j})\right|
\le\sum_{j=1}^{Q_n}H(\mathcal{A}^{\Delta_j})
\le c_4 Q_n^{\frac\alpha2+1}
\le c_5n^{\frac\alpha3+\frac23}
$$
as $H(\mathcal{A}^{\Delta_j})=\mathcal{O}(\Delta_j)=\mathcal{O}(n_j^\alpha)=\mathcal{O}(j^{\frac{\alpha}2})$.
The RHS is bounded by $n^{1-\delta_1}$ as $\delta_1\le\frac{1-\alpha}3$.
\qed

\vspace{2mm}

\noindent With $y_j=J_{n_j}\circ T^{N_j}$ we get $\mathbb{E}|y_j|^\ell\le n_j^\frac{\ell}2\le j^\frac{\ell}4$
by Proposition~\ref{higher.moments}. By Proposition~\ref{variance.convergence} then
$$
\sigma^2(y_j)=\mathbb{E}(y_j^2)=n_j\sigma^2+\mathcal{O}(n_j^{1-\eta})
$$
where $\eta=\frac{p(w-2)-5w+1}{2p(w-2)+2w+2}$. This then implies (as $n_j\sim\sqrt{j}$)
\begin{eqnarray}
\sum_{j=1}^{Q_n}\mathbb{E}(y_j^2)&=&\sigma^2\sum_jn_j+\sum_{j=1}^{Q_n}\mathcal{O}(n_j^{1-\eta})
\notag\\
&=&\sigma^2n+\mathcal{O}\left(\sum_{j=1}^{Q_n}n_j^\alpha\right)
+\sum_{j=1}^{Q_n}\mathcal{O}\left(j^{\frac{1-\eta}2}\right)\notag\\
&=&n\sigma^2+\mathcal{O}(n^{\frac{\alpha+2}3}+n^{1-\frac\eta3})\label{variances}
\end{eqnarray}
for $\delta\le\min\left(\frac{p(w-2)-5w+1}{2p(w-2)+2w+2},\frac{1-\alpha}3\right)$ as 
$n-\sum_jn_j=\sum_{j=1}^{Q_n}\mathcal{O}(j^\frac\alpha2)=\mathcal{O}(Q_n^{\frac\alpha2+1})
=\mathcal{O}(n^\frac{\alpha+2}3)$.

\begin{lemma}\label{sum.of.variances} If $\ell\ge5$ and $p\ge\frac{w(\ell+2)-\ell}{w-\ell}$. Then
$$
\sum_{j=1}^{Q_n}y_j^2=n\sigma^2+\mathcal{O}(n^{1-\delta_2})
$$
almost surely, for any $\delta_2<\min(\frac13,(1-\frac4\ell)\frac{\alpha p}6)$.
\end{lemma}

\noindent {\bf Proof.} To use Gal-Kuksma's estimate as given in Lemma~A1 of \cite{PS} directly we put 
$$
x_j=y_j^2-\mathbb{E}(y_j^2)
$$
and then obtain for $1\le m<m'\le Q_n$:
\begin{eqnarray*}
\mathbb{E}\left(\sum_{j=m}^{m'}x_j\right)^2
&=&\mathbb{E}\left(\sum_{j=m}^{m'}y_j^2-\mathbb{E}(y_j^2)\right)^2\\
&=&\sum_{j=m}^{m'}\left(\mathbb{E}(y_j^4)-(\mathbb{E}(y_j^2))^2\right)
+\sum_{j\not=i}\left(\mathbb{E}(y_j^2y_i^2)-\mathbb{E}(y_j^2)\mathbb{E}(y_j^2)\right)\\
&=&I+II.
\end{eqnarray*}
For the second term, $i\not=j$, we use Lemma~7.2.1 from~\cite{PS}: if $i<j$ then
$$
\left|\mathbb{E}(y_j^2y_i^2)-\mathbb{E}(y_j^2)\mathbb{E}(y_j^2)\right|
\le\beta(\Delta_j)^\frac1u\|y_i^2\|_r\|y_j^2\|_s
$$
where $\frac1u+\frac1r+\frac1s=1$. For the terms on the RHS we get by Proposition~\ref{higher.moments}
$$
\|y_i^2\|_s=\left(\mathbb{E}(y_i^{2s})\right)^\frac1s\le c_1(n_i^s)^\frac1s\le c_1\sqrt{i}
$$
under the assumption that $2s, 2r\le\ell$ (which requires $p\ge\frac{w(\ell+2)-\ell}{w-\ell}$) and obtain
$$
\left|\mathbb{E}(y_i^2y_j^2)-\mathbb{E}(y_j^2)\mathbb{E}(y_j^2)\right|
\le\beta(j^\frac\alpha2)^\frac1un_in_j=\mathcal{O}(j^{-\frac{\alpha p}2}\sqrt{ij})
$$
Thus for any $1\le m<m'\le Q_n$ we get
\begin{eqnarray*}
II&\le&c_2\sum_{m\le i<j\le m'}j^{-\frac{\alpha p}{2u}} \sqrt{ij}\\
&\le&c_2\sum_{i=m}^{m'-1}\sqrt i\sum_{j=i+1}^{m'}j^{\frac12-\frac12\frac{\alpha p}{u}}\\
&\le&c_3\sum_{i=m}^{m'-1}\sqrt i (m'^{\frac32-\frac12\frac{\alpha p}{u}}-i^{\frac32-\frac12\frac{\alpha p}{u}})\\
&\le&c_4(m'^\zeta-m^\zeta),
\end{eqnarray*}
where $\zeta=3-\frac12\frac{\alpha p}u$.
As $\mathbb{E}(y_j^4)=\mathcal{O}(n_j^2)$ we bound the first term~(I) using 
Proposition~\ref{higher.moments}
$$
I\le\sum_{j=m}^{m'}\left|\mathbb{E}(y_j^4)-(\mathbb{E}(y_j^2))^2\right|
=\sum_{j=m}^{m'}\mathcal{O}\left(n_j^2+n_j^2\right)
=\sum_{j=m}^{m'}\mathcal{O}(j)
=\mathcal{O}(m'^2-m^2).
$$
With $\zeta'=\max(\zeta,2)$ by~\cite{PS} Lemma~A.1 for any $\bar\delta>0$ there exist
a constant $c_5$ such that
$$
\sum_{j=1}^{Q_n}x_j
=\sum_{j=1}^{Q_n}\left(y_j^2-\mathbb{E}(y_j^2)\right)
\le c_5Q_n^\frac{\zeta'}2\log^{2+\bar\delta}Q_n
=\mathcal{O}(n^{1-\delta_2})
$$
for any $\delta_2<\min(\frac13,(1-\frac4\ell)\frac{\alpha p}6)$
almost surely where we have chosen $s=r=\frac{\ell}2$ which implies $\frac1u=1-\frac1r-\frac1s=1-\frac4\ell$
which is positive as $\ell\ge5$ by assumption.
Thus 
$$
\sum_{j=1}^{Q_n}y_j^2=\sum_{j=1}^{Q_n}\mathbb{E}(y_j^2)+\mathcal{O}(n^{1-\delta_2})
=\sigma^2n+\mathcal{O}(n^{1-\delta_2})
$$
almost surely.
\qed

\vspace{2mm}

\noindent
We now do the Martingale decomposition. Let $\mathscr{F}_j=\sigma(\mathcal{A}^{h_j})$,
where $h_j=N_j+n_j$. Then $\mathscr{F}_j\subset\mathscr{F}_{j+1}\subset\mathscr{F}_{j+2}\subset\cdots$
and introduce $Y_j$ by $y_j=Y_j+u_j-u_{j+1}$ where 
$$
u_j=\sum_{k=0}^\infty\mathbb{E}(y_{j+k}|\mathscr{F}_{j-1}).
$$

\begin{lemma}\label{coboundaries} There exists a constant $C_6$ such that 
$$
\|u_j\|_q^q\le C_6 j^{\frac14-\frac{\alpha p}2(1-\frac1q-\frac1\ell)}
$$
for $q<\frac14(3\alpha p(1-\frac1q-\frac1\ell)-1)$, provided $p\ge\frac{w(\ell+2)-\ell}{w-\ell}$.
\end{lemma}

\noindent {\bf Proof.} Since $\mathscr{F}_{j-1}$ `lives' on the first  $h_{j-1}$ coordinates 
 we obtain by~\cite{PS}~Lemma~7.2.1
\begin{eqnarray*}
\mathbb{E}\left(\left|\mathbb{E}(y_{j+k}|\mathscr{F}_{j-1})\right|^q\right)
&=&\mathbb{E}\left(|y_{j+k}|\cdot\left|\mathbb{E}(y_{j+k}|\mathscr{F}_{j-1})\right|^{q-1}\right)\\
&\le&\|y_{j+k}\|_r\left\|\left|\mathbb{E}(y_{j+k}|\mathscr{F}_{j-1})\right|^{q-1}\right\|_s
\beta^\frac1u(h_{j+k}-n_{j+k}-h_{j-1})
\end{eqnarray*}
for $\frac1r+\frac1s+\frac1u=1$, where
\begin{align*}
h_{j+k}-h_{j-1}&=\sum_{i=j}^{j+k-1}\left(n_i+\Delta_i\right)+\Delta_{j-1}\\
&=\mathcal{O}(1)\left(\sum_{i=j}^{j+k-1}\left(\sqrt i+i^\frac\alpha2\right)+(j-1)^\frac\alpha2\right)\\
&=\mathcal{O}(1)\left((j+k-1)^\frac32-j^\frac32+j^\frac\alpha2\right)
\end{align*}
Now 
$$ \left\|\left|\mathbb{E}(y_{j+k}|\mathscr{F}_{j-1})\right|^{q-1}\right\|_s
=\left(\mathbb{E}\left|\mathbb{E}(y_{j+k}|\mathscr{F}_{j-1})\right|^{s(q-1)}\right)^\frac1s
$$
and let $s$ be so that $s(q-1)=q$ i.e.\ $s=\frac{q}{q-1}$.
Then
$$
\left(\mathbb{E}\left|\mathbb{E}(y_{j+k}|\mathscr{F}_{j-1})\right|^{s(q-1)}\right)^\frac1q
\le\|y_{j+k}\|_r\beta^\frac1u(h_{j+k}-n_{j+k}-h_{j-1})
$$
where we used that $1-\frac1s=1-\frac{q-1}q=\frac1q$.
By Proposition~\ref{higher.moments} $\mathbb{E}|y_{j+k}|^r\le c_1n_{j+k}^\frac{r}2\le c_1(j+k)^\frac{r}4$ 
which implies $\|y_{j+k}\|_r\le c_2(j+k)^\frac14$ provided $r\le\ell$.
Since $(j+k-1)^\frac32-j^\frac32\ge k^\frac32$ one has
$$
\beta(h_{j+k}-n_{j+k}-h_{j-1})\le\beta\left(k^\frac32+j^\frac\alpha2\right)
\le c_3\left(k^\frac32+j^\frac\alpha2\right)^{-p}
\le c_3\min\left(j^{-\frac{\alpha p}2},k^{-\frac{3p}2}\right)
$$
and
$$
\left\|\mathbb{E}(y_{j+k}|\mathscr{F}_{j-1})\right\|_q^q
\le c_2c_3(j+k)^\frac14\left(k^\frac32+j^\frac\alpha2\right)^{-p}
$$
and consequently with $r=\ell$ (largest possible) and $\frac1u=1-\frac1q-\frac1\ell$
$$
\|u_j\|_q
\le c_2c_3\sum_{k=0}^\infty\left((j+k)^\frac14\left(k^\frac32+j^\frac\alpha2\right)^{-\frac{p}u}\right)^\frac1q
=\mathcal{O}( j^{\frac14-\frac{\alpha p}2(1-\frac1q-\frac1\ell)})^\frac1q
$$
provided $p$ is large enough, so that $\frac1{4q}(1-\frac{2\alpha p}u)<-1$ which is satisfied
by assumption.
\qed

\begin{lemma} \label{Martin.variance}
Let $\ell\ge5$ and assume $p\ge\max(\frac3{2\alpha},\frac9{\alpha(\ell-2)}\frac{w(\ell+2)-\ell}{w-\ell})$.
Then
$$
\sum_{j=1}^{Q_n}Y_j^2=\sigma^2n+\mathcal{O}(n^{1-\delta_3})
$$
almost surely for every $\delta_3<\min\left(\frac13,\frac23\alpha p-1,(1-\frac4\ell)\frac{\alpha p}6\right)$.
\end{lemma}

\noindent {\bf Proof.}
Put $v_j=u_j-u_{j+1}$, then $Y_j=y_j-v_j$ and $Y_j^2=y_j^2-2v_jy_j+v_j^2$.
We estimate the two terms separately. 
First by Lemma~\ref{coboundaries} with ($q=2$)
$$
\mathbb{E}(v_j^2)=\mathbb{E}(u_j-u_{j+1})^2
\le 4(\max(\|u_j\|_2,\|u_{j+1}\|_2))^2
\le c_1j^{\frac12-\alpha p}
$$
($c_1\le4C_9$) provided $p(\frac12-\frac1\ell)>\frac9{2\alpha}$ (which is satisfied by assumption) and, 
using Cauchy-Schwarz,
$$
\mathbb{E}|y_jv_j|\le\|y_j\|_2\|v_j\|_2
\le c_1\sqrt{\mathbb{E}(y_j^2)}j^{\frac14-\frac12\alpha p}
\le c_2n_j^\frac12j^{\frac14-\frac12\alpha p}
\le c_2j^{\frac12-\frac12\alpha p}.
$$
Then
$$
\mathbb{E}\left(\sum_{j=1}^{Q_n}v_j^2\right)=\sum_{j=1}^{Q_n}\mathbb{E}(v_j^2)
\le c_1\sum_{j+1}^{Q_n}j^{\frac12-\alpha p}
\le c_4Q_n^{\frac32-\alpha p}
\le c_5n^{1-\frac23\alpha p}
$$
and with $\epsilon_n=n^{1-\delta_3'}$
$$
\mathbb{P}\left(\sum_{j=1}^{Q_n}v_j^2\ge\epsilon_n\right)
\le\frac1{\epsilon_n}\mathbb{E}\left(\sum_{j=1}^{Q_n}v_j^2\right)
\le c_5\frac{n^{1-\frac23\alpha p}}{\epsilon_n}
\le c_5n^{\delta_3'-\frac23\alpha p}.
$$
If $\delta_3'-\frac23\alpha p<-1$ then by Borel-Cantelli
$$
\sum_{j=1}^{Q_n}v_j^2<n^{1-\delta_3'}
$$
for all large enough $n$ almost surely. It now follows from Lemma~\ref{sum.of.variances} that
$$
\sum_{j=1}^{Q_n}Y_j^2=\sum_{j=1}^{Q_n}y_j^2+\mathcal{O}(n^{1-\delta_3'})
=\sigma^2n+\mathcal{O}(n^{1-\delta_3})
$$
almost surely for all $\delta_3<\min\left(\frac13,\frac23\alpha p-1,(1-\frac4\ell)\frac{\alpha p}6\right)$.
\qed

\begin{lemma} \label{Martin.conditional} Let $\ell\ge3$, then
$$
\sum_{j=1}^{Q_n}\left(\mathbb{E}(Y_j^2|\mathscr{F}_{j-1})-Y_j^2\right)
\le n^{1-\delta_4}
$$
almost surely, for any $\delta_4<\frac23(1-\frac2\ell)$.
\end{lemma}

\noindent {\bf Proof.} If we put $R_j=\mathbb{E}(Y_j^2|\mathscr{F}_{j-1})-Y_j^2$, a 
Martingale difference, then by Minkowsky and Proposition~\ref{higher.moments}
$$
\mathbb{E}|R_j^q|\le\mathbb{E}|Y_j|^{2q}
\le \mathbb{E}|y_j|^{2q}+\mathbb{E}|v_j|^{2q}
\le C_5n_j^q+C_6j^{(\frac14-\frac12\alpha p)q}
\le c_1 j^\frac{q}2
$$
provided $2q\le\ell$, and therefore
$$
\sum_{j=1}^\infty j^{-q\gamma}\mathbb{E}|R_j|^q
\le c_1\sum_{j=1}^\infty j^{-q\gamma}j^\frac{q}2<\infty
$$
if $q(\frac12-\gamma)<-1$ i.e.\ $\gamma>\frac12+\frac1q$.
Then $\sum_j j^{-\gamma}R_j$
converges almost surely and therefore
by Kronecker's lemma
$$
\sum_{j=1}^{Q_n}R_j=\mathcal{O}(Q_n^\gamma)=\mathcal{O}(n^{\frac23\gamma})
$$
almost surely if $\frac23\gamma=1-\delta_4$, where $\delta_4<\frac23(1-\frac1q)$
and where we chose $q=\frac\ell2$ (the largest possible value). (Note $\frac23(1-\frac1q)>\frac13$.)
\qed

\vspace{3mm}

\noindent {\bf Proof of Theorem~\ref{ASIP}. }
By the Skorokhod representation theorem there exist $T_j$ such 
$$
\sum_{j=1}^{Q_n}Y_j=X\left(\sum_{j=1}^{Q_n}T_j\right)
$$
almost surely, where $X$ is the Brownian motion and where 
$\mathbb{E}(T_j|\mathscr{F}_{j-1})=\mathbb{E}(Y_j^2|\mathscr{F}_{j-1})$ a.s.\
and $\mathbb{E}(T_j^2)\le\mathbb{E}|Y_j|^{2q}$ for $q>1$.
Then, we conclude as in Philip and Stout~\cite{PS}
$$
\sum_{j=1}^{Q_n}T_j-\sigma^2n
=\sum_{j=1}^{Q_n}\left(T_j-\mathbb{E}(T_j|\mathscr{F}_{j-1})\right)
+\sum_{j=1}^{Q_n}\left(\mathbb{E}(T_j|\mathscr{F}_{j-1})-Y_j^2\right)
+\sum_{j=1}^{Q_n}Y_j^2-\sigma^2n.
$$
For the first term on the RHS we use that 
$$
\mathbb{E}\left|T_j-\mathbb{E}(T_j|\mathscr{F}_{j-1})\right|^2\le\mathbb{E}|T_j|^q\le \mathbb{E}|Y_j|^{2q},
$$
for the second term we use Lemma~\ref{Martin.conditional} and the third term was estimated in 
Lemma~\ref{Martin.variance}. Notice that since $\ell\ge5$ we get that $\frac23(1-\frac2\ell)>\frac13$
(Lemma~\ref{Martin.conditional}). Hence 
$$
\delta<\min\left(\frac{p(w-2)-5w+1}{2p(w-2)+2w+2},
\sup_{\alpha\in(0,1)}\min \left(\frac{1-\alpha}3,\frac{\alpha p}{30},\frac{2\alpha p}3-1\right) \right).
$$
To get the statement of the theorem let us look at the second term and notice that the last 
two entries in it agree when $\alpha p=\frac{30}{19}$ and there they produce the value $\frac1{19}$.
Hence the supremum is realised at $\alpha=\frac{10}{p+10}$ and its value equals
$\frac13-\frac{10}{3p+30}$. Consequently
\begin{equation}\label{delta.bound}
\delta<\min\left(\frac{p(w-2)-5w+1}{2p(w-2)+2w+2},\frac13-\frac{10}{3p+30} \right).
\end{equation}
In the first term we can for instance put $\ell=5$ and $p=7$ and thus obtain the value $\frac18$
which is smaller than the second term.
Therefore we can use any $\delta\le\frac18$.
\qed

\vspace{3mm}

\noindent {\bf Proof of Corollary~\ref{corollary}.} In part~(i) if we let $w\to\infty$ then 
$\delta<\min\left(\frac12-\frac3{p+1},\frac13-\frac{10}{3p+30} \right)$.
For part~(ii) we let in~\eqref{delta.bound} $p$ go to infinity which leads to  the condition
$\delta<\min(\frac12,\frac13)$.
\qed



\begin{thebibliography}{99}

\bibitem{Aba} M Abadi: Exponential Approximation for
Hitting Times in Mixing Stochastic Processes; Mathematical Physics Electronic
Journal~{\bf 7} (2001).



\bibitem{Bressaud} X Bressaud: Subshifts on an infinite alphabet; Ergod.\ Th.\
\& Dynam.\ Sys.\ ~{\bf 19} (1999), 1175-1200.

\bibitem{BK} M Brin and A Katok: On local entropy; {\it Geometric Dynamics} 30--38,
Springer LNM~\#1007, 1983.

\bibitem{Bro} A Broise: {\it Transformations dilatantes de l'intervalle et
th\'eor\`emes limites}; Asterisque~\#238, 1996.

\bibitem{BV} H Bruin and S Vaienti: Return times for unimodal maps;
submitted to Forum Math.\


\bibitem{ch} N Chernov: Limit theorems and Markov approximations for
chaotic dynamical systems; Prob.\ Th.\ Rel.\ Fields~{\bf 101} (1995), 321--362.



\bibitem{CGS} P Collet, A Galves and B Schmitt: Fluctuations of
repetition times for Gibbsian sources; Nonlinearity~{\bf 12} (1999), 1225--1237.

\bibitem{Dou} P Doukhan: {\it Mixing: Properties and examples}; Lecture Notes
in Statistics~85, Springer 1995.

\bibitem{FHV} P Ferrero, N Haydn and S Vaienti: Entropy fluctuations
for parabolic maps; Nonlinearity~{\bf 16} (2003), 1203--1218.


\bibitem{GS} A Galves and B Schmitt: Inequalities for hitting
times in mixing dynamical systems; Random and Computational Dynamics~{\bf 5} (1997), 337--348.

\bibitem{Gor} M Gordin: The central limit theorem for stationary processes;
Soviet Math.\ Doklady~{\bf 10} (1969), 1174--1176.

\bibitem{Han} Guanyue Han: Limit Theorems in Hidden Markov Models;
IEEE Trans.\ Inform.\ Theory {\bf 59(3)} (2013), 1311--1328.

\bibitem{Hay} N Haydn: The Central Limit Theorem for uniformly strong mixing measures;
Stochastics and Dynamics 2012.

\bibitem{HV2} N Haydn and S Vaienti: Fluctuations of the metric entropy 
for mixing measures; Stochastics and Dynamics~4 (2004), 595--627

\bibitem{HV3} N Haydn and S Vaienti: The distribution of the measure of 
cylinders for non-Gibbsian measures; {\it Complex Dynamics and
Related Topics} 147--162, New Studies in Advanced Mathematics~\#5, 2003.

\bibitem{Ibr} I A Ibragimov: Some limit theorems for stationary processes;
 Theory Prob.\ Appl.~{\bf 7} (1962), 349--382.

\bibitem{Kon} I Kontoyiannis: Asymptotic Recurrence and Waiting
Times for Stationary Processes; J. Theor.\ Prob.~{\bf 11} (1998), 795--811.

\bibitem{Liv} C Liverani: Central Limit theorem for deterministic systems;
Intern.\ Congress on Dyn.\ Syst., Montevideo 1995 (Proc.\ Research Notes in
Math.\ Series), Pitman (1996), 56--75.

\bibitem{Man} R Ma\~n\'e: {\it Ergodic Theory and Differentiable Dynamics},
Springer 1987.


\bibitem{NW} A Nobel and A Wyner: A Recurrence Theorem for Dependent 
Processes with Applications to Data Compression; IEEE Vol.~{\bf 38(5)} (1992), 1561--1564.

\bibitem{OW} Ornstein and Weiss; Entropy and Data Compression Schemes;
IEEE Transactions on Information Theory~{\bf 39} (1993), 78--83.

\bibitem{OW2} Ornstein and Weiss; Entropy and Recurrence Rates for Stationary
Random Fields; IEEE Transactions on Information Theory~{\bf48(6)} (2002), 1694--97.

\bibitem{Pac} F Paccaut: {\it Propri\'et\'es Statistiques de Syst\`emes
Dynamiques Non Markovian}; Th\`ese (Doctorat) Dijon 2000.

\bibitem{Pen} F P\`ene: Rates of Convergence in the CLT for
Two-Dimensional Dispersive Billiards; Commun.\ Math.\ Phys.~{\bf 225} (2002), 91--119.

\bibitem{Petersen} K Petersen: {\it Ergodic Theory}: Cambridge studies in advanced
mathematics~\#2, 1983.


\bibitem{PS} W Philipp and W Stout: {\it Almost sure invariance principles for 
partial sums of weakly dependent random variables}; AMS Memoirs vol.~2 
issue~2 No.~161, 1975.



\bibitem{Rom} V I Romanovskii: {\it Discrete Markov Chains}; Wolters-Noordhoff
Publishing Groningen 1970.

\bibitem{R1} M Rosenblatt: A central limit theorem and a strong mixing
condition; Proc.\ Nat.\ Acad.\ Sci.\ USA~{\bf 42} (1956), 43--47.

\bibitem{R2} M Rosenblatt: {\it Markov Processes. Structure and Asymptotic
Behavior}; Springer, Grundlagen~\#184, 1971.

\bibitem{Sarig} O Sarig: Thermodynamic formalism for countable Markov shifts;
Ergod.\ Th.\ \& Dynam.\  Syst.~{\bf 19} (1999), 1565--1593.

\bibitem{Wal} P Walters: {\it An Introduction to Ergodic Theory}; Springer-Verlag 1981.

\bibitem{WZ} A Wyner and J Ziv: Some asymptotic properties of the entropy of a stationary ergodic
 data source with applications to data compression; IEEE Trans.\ Inf.\ Theory {\bf 35} (1989), 1250--8.


\bibitem{Yus} A Yushkevich: On limit theorems connected with the concept
of the entropy of Markov chains (in Russian); Uspehi Math.\ Nauk~{\bf 8} (1953) 177--180.


\end{thebibliography}
\end{document}